\magnification=1090
\hsize=12cm
\hoffset=1,2cm
\voffset=5mm
\baselineskip=13pt
\overfullrule=0mm 
\font\a=msbm10
\font\b=cmr17 at 18pt

\font\f=cmcsc10
\font\ff=cmcsc10 at 8pt
\font\bbb=cmr10 at 8pt
\font\bbbb=cmr17 at 12pt


\def\build#1^#2{\mathrel{\mathop{\kern 0pt#1}\limits^{#2}}}
\centerline{\b Homotopical variations and high-dimensional}
\medskip
\centerline{\b Zariski-van Kampen theorems}
\vskip 6mm
\centerline{\bbbb D. CH\'ENIOT $^{(1)}$ and C. EYRAL
$^{(2)}$}
\vskip 6mm

\baselineskip=11pt
\centerline{\bbb \item{(1)} Laboratoire d'Analyse, Topologie et Probabilit\'es
(UMR CNRS 6632)} 
\centerline{\bbb Centre de Math\'ematiques et d'Informatique, Universit\'e de
Provence}
\centerline{\bbb 39 rue Joliot-Curie, 13453 Marseille C\'edex 13, France}
\centerline{\bbb E-mail: cheniot@gyptis.univ-mrs.fr}\vskip 3mm  

\centerline{\bbb \item{(2)} Department of Mathematics}
\centerline{\bbb The Abdus Salam International Centre for Theoretical Physics}
\centerline{\bbb Strada Costiera 11, 34014 Trieste, Italy}
\centerline{\bbb E-mails: eyralc@ictp.trieste.it or eyralchr@yahoo.com} 
\vskip 6mm

\baselineskip=13pt
\centerline{\bbbb 11 April, 2002}\vskip 8mm

\leftskip=1,5cm\rightskip=1,5cm
\baselineskip=11pt
\centerline{\ff Abstract}
{\bbb In 1933, van Kampen described the fundamental groups of
the complements of plane complex projective algebraic curves.
Recently, Ch\'eniot-Libgober proved an analogue of this result
for higher homotopy groups of the complements of complex
projective hypersurfaces with isolated singularities. Their
description is in terms of some ``homotopical variation
operators''. We generalize here the notion of ``homotopical
variation'' to (singular) quasi-projective varieties. This is
a first step for further generalizations of van Kampen's theorem.
A conjecture, with a first approach, is stated in the special case of
non-singular quasi-projective varieties.}

\leftskip=0cm\rightskip=0cm
\vskip3mm

{\bbb \noindent Mathematics Subject Classification 2000: 14F35
(Primary), 14D05, 32S50, 55Q99 (Secondary).

 }
\vskip 1cm
\baselineskip=13pt

\centerline{\f Introduction}\vskip 3mm

Our topic is best understood in the general frame of
Lefschetz type theorems. Let $X:=Y\setminus Z$, where $Y$~is an
algebraic subset of complex projective space~$\hbox{\a P}^n$,
with $n\geq 2$, and $Z$~an algebraic subset of~$Y$ (such
an~$X$ is called an (embedded) quasi-projective variety). Let
${\cal L}$~be a projective hyperplane of~$\hbox{\a P}^n$. 
The non-singular quasi-projective version of the Lefschetz
Hyperplane Section Theorem, proved by Hamm-L\^e and
Goresky-MacPherson ({\it cf.}~[HL1] and~[GM1,2]), asserts that
if $X$~is non-singular and $\cal L$~generic, then the natural
maps (between homology and homotopy groups, respectively)
$$H_q({\cal L}\cap X)\rightarrow H_q(X)\quad\hbox{and}\quad
\pi_q({\cal L}\cap X,*)\rightarrow \pi_q(X,*)$$
are bijective for $0\leq q\leq d-2$ and surjective for
$q=d-1$, where $d$~is the least (complex) dimension of the
irreducible components of~$Y$ not contained in~$Z$. In the
special case $Y=\hbox{\a P}^n$, that is for the complement
$X=\hbox{\a P}^n\setminus Z$ of a projective variety, the
bounds can be improved by~$c-1$, where $c$~is the least
codimension of the irreducible components of~$Z$ ({\it
cf}.~[C2]): then the above maps are bijective for $0\leq q\leq
n+c-3$ and surjective for $q=n+c-2$ (there is no improvement
if~$c=1$).

The question now arises of determining the kernel of these
maps in dimension~$d-1$ (resp.~$n+c-2$ for a complement). For
this purpose, it is classical (at least when $Z=\emptyset$) to
consider~${\cal L}$ as a member of a pencil~${\cal P}$ of
hyperplanes of~$\hbox{\a P}^n$ with axis a generic
$(n-2)$-plane~${\cal M}$. The sections of~$X$ by all the
hyperplanes of such a pencil are isotopic to the one
by~${\cal L}$ with the exception of the sections by a finite
number of exceptional hyperplanes~$({\cal L}_i)_i$. For each~$i$,
there are some homomorphisms
$$\hbox{var}_{i,q}\colon\ H_q({\cal L}\cap X,{\cal M}\cap X)
\rightarrow H_{q}({\cal L}\cap X),\hbox{ for all $q$},$$
called ``homological variation operators'', defined by
patching each relative cycle on ${\cal L}\cap X$ modulo
${\cal M}\cap X$ with its transform by monodromy
around~${\cal L}_i$ ({\it cf.}~[C3]). These are defined even
if $X$~is singular. In {\it loc.~cit.}, the first author showed
that if $X$~is non-singular, then the kernel of the natural
epimorphism
$$H_{d-1}({\cal L}\cap X)\rightarrow H_{d-1}(X)\leqno(E)$$
is equal to the sum of the images of the homological variation
operators\break $(\hbox{var}_{i,d-1})_i$ associated to the
exceptional members~$({\cal L}_i)_i$ of the pencil. In the
case of a complement $X=\hbox{\a P}^n\setminus Z$, the same is
true for the epimorphism
$$H_{n+c-2}({\cal L}\cap(\hbox{\a P}^n\setminus Z))\rightarrow
H_{n+c-2}(\hbox{\a P}^n\setminus Z)\leqno(E')$$
with $n+c-2$ in place of $d-1$.
Combined with the Hyperplane Section Theorem quoted above,
this gives a natural isomorphism
$$H_{d-1}({\cal L}\cap X)\Big/\sum_i\hbox{Im var}_{i,d-1}
\buildrel\sim\over\longrightarrow H_{d-1}(X)\leqno (1)$$
and, in the case of a complement,
$$H_{n+c-2}({\cal L}\cap(\hbox{\a P}^n\setminus
Z))\Big/\sum_i\hbox{Im var}_{i,n+c-2}
\buildrel\sim\over\longrightarrow
H_{n+c-2}(\hbox{\a P}^n\setminus Z)\leqno (1')$$
(which coincides with a special case of~(1) when~$c=1$).

This is different from the classical point of view of the
Second Lefschetz Theorem ({\it cf.}~[Lef], [Wa] and~[AF]) where
the kernel is expressed in terms of ``vanishing cycles'', that
is $(d-1)$-cycles of~${\cal L}\cap X$ which vanish when ${\cal
L}$~tends to some~${\cal L}_i$. But the classical theorem
applies only to the case where $X$~is non-singular {\it
closed\/} (i.e.,
$Z=\emptyset$) while formula~(1) applies to both the closed
and non-closed cases, provided $X$~is non-singular. Thus
isomorphisms~(1) and~$(1')$ are generalizations of the Second
Lefschetz Theorem.

But isomorphisms (1) and~$(1')$ may also be considered as
generalized homological forms of the classical
Zariski-van~Kampen theorem on curves. Recall that this theorem
expresses the fundamental group of the complement of an
algebraic curve in~$\hbox{\a P}^2$ by generators and
relations. The generators are loops in a generic line, around
its intersection points with the curve. The relations are
obtained by considering a generic pencil containing this line:
each loop must be identified with its transforms by monodromy
around the exceptional members of the pencil ({\it cf.}~[Za],
[vK] and~[C1]).

The first true (homotopical) generalization of van~Kampen's
theorem to higher dimensions was given by Libgober ({\it
cf}.~[Li]). It applies to the $(n-1)$-st homotopy group of the
complement of a hypersurface with isolated singularities
in~$\hbox{\a C}^n$ behaving well at infinity. In
this case, if $n\geq3$, the fundamental group is~$\hbox{\a Z}$
and the $(n-1)$-st homotopy group is the first higher
non-trivial one as explained in~[Li], and at
the same time the first not preserved by hyperplane section.
Libgober also showed how the projective case can be deduced
from the affine case. The projective version of Libgober's
theorem can be stated in the previous frame as follows.
With the notation above, assume that~$Y=\hbox{\a P}^n$ with
$n\geq 3$, and $Z=H$ where~$H$ is an algebraic hypersurface of
$\hbox{\a P}^n$ having only isolated singularities. Consider a
generic pencil of hyperplanes of~$\hbox{\a P}^n$ as above. Let
$*\in {\cal M}\cap(\hbox{\a P}^n\setminus H)$ be a base point
in the axis of the pencil. For each~$i$, we denote by  
$$h_{i,q}\colon\ \pi_q({\cal L}\cap (\hbox{\a
P}^n\setminus H),*)\rightarrow\pi_q({\cal L}\cap(\hbox{\a
P}^n\setminus H),*),\hbox{ for all }q\geq 1,$$
the isomorphism induced by a geometrical monodromy
around~${\cal L}_i$ leaving fixed the points of
${\cal M}\cap(\hbox{\a P}^n\setminus H)$. Then there is a
natural isomorphism
$$\eqalign{\pi_{n-1}({\cal L}\cap (\hbox{\a P}^n\setminus
&H),*)\Big/\Bigl(\hbox{Im }(h_{1,n-1}-\hbox{id}),\hbox{Im
}{\cal D}_{1,n-1},\cr
&\ldots,\hbox{Im }(h_{N,n-1}-\hbox{id}),\hbox{Im
}{\cal D}_{N,n-1}\Bigr)
\buildrel\sim\over\longrightarrow\pi_{n-1} (\hbox{\a
P}^n\setminus H,*),\cr}\leqno(2)$$
where $N$~is the number of the exceptional hyperplanes and where
each
$${\cal D}_{i,n-1}\colon\ \pi_{n-2}({\cal L}_i\cap (\hbox{\a
P}^n\setminus H),*)\rightarrow \pi_{n-1}({\cal L}\cap
(\hbox{\a P}^n\setminus H),*)\big / 
\hbox{Im }(h_{i,n-1}-\hbox{id})$$
is some homomorphism called ``degeneration operator''.

Later in~[CL], Ch\'eniot and Libgober showed that the
combination (ordinary variations by monodromies - degeneration
operators) of~[Li] is related to some homotopical variation
$${\cal V\hskip -0,8mm AR}_{i,n-1}\colon\ \pi_{n-1}({\cal
L}\cap (\hbox{\a P}^n\setminus H), {\cal M}\cap (\hbox{\a
P}^n\setminus H),*)\rightarrow\pi_{n-1}({\cal L}\cap (\hbox{\a
P}^n\setminus H),*)$$
defined from the homological variation $\hbox{var}_{i,n-1}$
of~[C3] acting in a pencil supported by the universal
covering of $\hbox{\a P}^n\setminus H$. In particular,
(2)~yields a natural isomorphism
$$\pi_{n-1}({\cal L}\cap
(\hbox{\a P}^n\setminus H),*)\Big/\sum_i
\hbox{Im }{\cal V\hskip -0,8mm
AR}_{i,n-1}\buildrel\sim\over\longrightarrow
\pi_{n-1}(\hbox{\a P}^n\setminus H,*).\leqno (3)$$
Isomorphism~(3) can also be deduced from~(1) applied to the
universal covering of $\hbox{\a P}^n\setminus H$ ({\it
cf}.~[CL]). 

Isomorphism~(3) provides a homotopy analogue of
isomorphisms~(1) and~$(1')$ in the special case of complements of
projective hypersurfaces with isolated singularities. Thus,
this high-dimensional Zariski-van~Kampen theorem is also a
homotopy version of the Second Lefschetz Theorem for this
special case. But the definition of homotopical variation
operators in~[CL] (as well as the definition of degeneration
operators in~[Li]) relies heavily upon the special topology of
complements of hypersurfaces with isolated singularities and
does not make sense in a more general setting. The need then
arises for an alternative definition which could lead to a more
general homotopy analogue of isomorphisms~(1) and~$(1')$. It
must be said that in a more general situation, the considered
homotopy group may not be the first higher non-trivial one but
it remains the first not preserved by generic hyperplane
section.  The aim of this article is to give such an alternative
definition and to state a reasonable conjecture generalizing~(3)
(and~(2)) with a first approach toward its proof.

We introduce here new homotopical variation operators,
$$\hbox{VAR}_{i,q}\colon\ \pi_{q}({\cal L}\cap X,{\cal M}\cap
X,*)\rightarrow \pi_{q}({\cal L}\cap X,*),\hbox{ for all }
q\geq 1,$$
which extend those of Ch\'eniot-Libgober and have several
advantages with respect to them. Our definition is purely
homotopical, that is, it does not go through homology. This
frees it from the special topology of the case considered by
Ch\'eniot and Libgober, which was precisely required to express
homotopy groups with the help of homology groups. In fact our
definition is valid for any (singular) quasi-projective variety
$X:=Y\setminus Z$. Our operators coincide with those of
Ch\'eniot-Libgober when the latter are defined. More
precisely, we show that, when $X:=\hbox{\a P}^n\setminus H$ and
$q=n-1$, with $n\geq 3$, and $H$~is an algebraic hypersurface
having only isolated singularities, then
$\hbox{VAR}_{i,n-1}={\cal V\hskip -0,8mm AR}_{i,n-1}$. In
particular, isomorphism~(3) is valid with $\hbox{VAR}_{i,n-1}$
in place of~${\cal V\hskip -0,8mm AR}_{i,n-1}$. We also show
that our homotopical operators are linked by Hurewicz
homomorphisms with the homological operators
$\hbox{var}_{i,q}$ of~[C3] and that they are equivariant under
the action of the fundamental group.

It is then natural to ask whether isomorphism~(1) when $X$~is
non-singular and isomorphism~$(1')$ are true for homotopy groups,
with our homotopical variation operators instead of the
homological ones. We show easily that the kernel of the homotopy
analogue of epimorphism~$(E)$ (resp.~$(E')$) contains the images of
operators $\hbox{VAR}_{i,d-1}$ (resp.~$\hbox{VAR}_{i,n+c-2}$).
Thus, there are well defined epimorphisms
$$\eqalign{
\pi_{d-1}({\cal L}\cap X,*)\Big/\sum_i\hbox{Im VAR}_{i,d-1}
 &\rightarrow\pi_{d-1}(X,*)\quad\hbox{if $d\geq3$},\cr
\pi_1({\cal L}\cap X,*)\bigg/\overline{\bigcup_i\hbox{Im
VAR}_{i,1}}
 &\rightarrow\pi_1(X,*)\quad\hbox{if $d=2$},\cr
}\leqno(4)
$$
where $\overline{\bigcup_i\hbox{Im VAR}_{i,1}}$ denotes the
normal subgroup generated by
$\bigcup_i\hbox{Im VAR}_{i,1}$, and a well-defined epimorphism
$$\pi_{n+c-2}({\cal L}\cap (\hbox{\a P}^n\setminus
Z),*)\Big/\sum_i\hbox{Im VAR}_{i,n+c-2}
 \rightarrow\pi_{n+c-2}(\hbox{\a P}^n\setminus Z,*)
\leqno(4')$$
when $n+c\geq4$ (if $c=1$, $(4')$~is only a special case of~(4)).
The question is whether these are isomorphisms. Notice that
when $n\geq3$ and $X=\hbox{\a P}^n\setminus H$ is the complement of a
hypersurface~$H$ with isolated singularities, (4)~and $(4')$~coincide
with isomorphism~(3) since our operators extend those of
Ch\'eniot-Libgober. Also, when $n=2$ and $X:=\hbox{\a P}^2\setminus C$
is the complement of a plane curve, the second row
of~(4)~coincides with the map that van~Kampen's
theorem asserts to be an isomorphism, as we shall see in Section~4.

We conjecture that the answer to our question is almost always yes:
\vskip 3mm

{\f Conjecture.} -- {\it Epimorphism\/~$(4)$, where $X$~is
non-singular, and epimorphism\/~$(4')$, unless $Z=\emptyset$, are
isomorphisms.}
\vskip3mm

We are far from having a proof of this conjecture. We shall only
give a very small step in its direction in the last section.

A detailed exposition from the origins to nowadays of the
questions mentioned in this introduction, like the Lefschetz
Hyperplane Section Theorem, the Second Lefschetz Theorem or
the van Kampen Theorem on curves, can be found in~[E3].

The content of this article is as follows. Section~1 is
devoted to some basic facts on generic pencils and
monodromies. In Section~2, we recall the definition of the
homological variation operators $\hbox{var}_{i,q}$ of~[C3]
which are used to define the homotopical variation operators
${\cal V\hskip -0,8mm AR}_{i,n-1}$ of~[CL]. The latter will be
described in Section~3. In Sections~4 and 5, we introduce the
generalized homotopical variation
operators~$\hbox{VAR}_{i,q}$, we give their elementary
properties and we prove that they coincide with the
Ch\'eniot-Libgober operators when the latter are defined. We
also prove that there are epimorphisms~(4) and~$(4')$ as stated
above. Finally, in Section~6, we give a first approach to the
conjecture we have formulated.
\medskip

{\f Notation 0.1.}~-- Throughout the paper, homology groups are
singular homology groups with integer coefficients, unless there
is an explicit statement of the contrary. We shall note the
homology class in a space $A$ of an (absolute) cycle $z$
by~$[z]_A$ and the homology class in $A$ modulo a subspace~$B$
of a relative cycle $z'$ by~$[z']_{A,B}$. If there is no
ambiguity, we shall omit the subscripts. If
$(A,B)$ is a pointed pair with base point $*\in B$, we shall
denote by~$F^q(A,B,*)$ the set of relative homotopy $q$-cells of
$A$ modulo $B$ based at $*$. These are maps from the
$q$-cube~$I^q$ to~$A$ with the face $x_q=0$ sent into~$B$ and all
other faces sent to~$*$ (as in~[St, \S15]). We designate by
$F^q(A,*)$ the set of absolute homotopy $q$-cells of $A$ based
at~$*$, that is maps from~$I^q$ to~$A$ sending the
boundary~$\dot I^q$ of~$I^q$ to~$*$. Given
$f\in F^q(A,B,*)$ (resp.~$F^q(A,*)$), the homotopy class of~$f$ 
in $A$ modulo
$B$ based at $*$ (resp.~in~$A$ based at $*$) will be denoted by
$\langle f\rangle_{A,B,*}$ (resp.~$\langle f\rangle_{A,*}$).
Again, if there is no ambiguity, we shall omit the subscripts.

\vskip 1cm

\centerline{\f 1. Generic pencils and monodromies}\vskip 3mm

Let~$X:=Y\setminus Z$, where $Y$ is a non-empty closed algebraic
subset of
$\hbox{\a P}^n$, with $n\geq 2$, and~$Z$ a proper closed
algebraic subset of
$Y$. Take a Whitney stratification~${\cal S}$ of $Y$ such
that~$Z$ is a union of strata ({\it cf.}~[Wh], [LT]), and
consider a projective hyperplane ${\cal L}$ of $\hbox{\a P}^n$
transverse to (the strata of) ${\cal S}$ (the choice of such a
hyperplane is generic). Denote by~$d$ the least dimension of the
irreducible components of~$Y$ not contained in
$Z$.

Notice that, when $X$ is non-singular, the application of the Lefschetz
Hyperplane Section Theorems mentioned in the introduction
is valid for ${\cal L}$. Indeed, the genericity of the hyperplane
required for these theorems can be specified as 
its transversality to a Whitney
stratification of~$Z$ ({\it cf}.~[HL2, Appendix], [GM2, end of the
proof of~II.5.1], [C2, Corollaire~1.2]).

Now, consider a pencil ${\cal P}$ of hyperplanes of
$\hbox{\a P}^n$ having~${\cal L}$ as a member and the axis~${\cal M}$
of which is transverse to~${\cal S}$ (the choice of
such an axis is generic inside ${\cal L}$). All the members of
${\cal P}$ are transverse to ${\cal S}$ with the exception of a
finite number of them~$({\cal L}_i)_i$, called {\it exceptional
hyperplanes}, for which, nevertheless, there are only a finite
number of points of non-transversality, all of them situated
outside of~${\cal M}$~({\it cf.}~[C2]). If necessary, one may
take the liberty of considering some ordinary members
of~${\cal P}$, different from ${\cal L}$, as exceptional ones. 

We parametrize the elements of~${\cal P}$ by the complex
projective line~$\hbox{\a P}^1$ as usual. Let $\lambda$ be the
parameter of~${\cal L}$ and, for each $i$, let $\lambda_i$ be the
parameter of~${\cal L}_i$. For each~$i$, take a small closed semi-analytic
disk~$D_i\subset\hbox{\a P}^1$ with centre~$\lambda_i$ together with
a point~$\gamma_i$ on its boundary. Choose the
$D_i$ mutually disjoint. Take also the image $\Gamma_i$ of a
simple real-analytic arc in $\hbox{\a P}^1$ joining $\lambda$
to~$\gamma_i$ and such that:
$\hbox{(i) }\Gamma_i\cap D_i=\{\gamma_i\}$;
$\hbox{(ii) }\Gamma_i\cap \Gamma_{i'}=\{\lambda\}\hbox{ if
}i\not=i'$;
$\hbox{(iii) }\Gamma_i\cap D_{i'}=\emptyset \hbox{ if
}i\not=i'$. Then, set
$$K_i:=\Gamma_i\cup D_i.$$
Finally, consider a loop~$\omega_i$ in the boundary $\partial
K_i$ of~$K_i$ starting from~$\lambda$, running along $\Gamma_i$
up to $\gamma_i$, going once counter-clockwise around the
boundary of~$D_i$ and coming along $\Gamma_i$ back to~$\lambda$.
\vskip 3mm

{\f Notation 1.1.}~-- For any subsets $G\subset \hbox{\a P}^n$
and $E\subset\hbox{\a P}^1$, note
$$G_E:=\bigcup_{\mu\in E} G\cap {\cal P}(\mu),$$
where ${\cal P}(\mu)$ is the member of ${\cal P}$ with parameter
$\mu$.
\vskip 3mm\goodbreak
The monodromies around each~${\cal L}_i$, more precisely above
each~$\omega_i$, proceed from the following lemma.
\vskip3mm

{\f Lemma 1.2} ({\it cf.}~[C3, Lemma 4.1]).~-- {\it For each $i$,
there is an isotopy
$$H\colon\ ({\cal L}\cap X)\times
I\rightarrow X_{\partial K_i}$$
such
that:}\smallskip \item{ (i)} $H(x,0)=x$, {\it for every $x\in {\cal
L}\cap X$};\smallskip  \item{(ii)} $H(x,t)\in X_{\omega_i(t)}$ {\it
for every} $x\in {\cal L}\cap X$ {\it and every} $t\in
I$;\smallskip \item{(iii)} {\it for every $t\in I$, the map} ${\cal
L}\cap X\rightarrow X_{\omega_i(t)},$ {\it  defined by
$x\mapsto H(x,t)$, is a homeomorphism;}
\smallskip   \item{(iv)} $H(x,t)=x$, {\it for every $x\in {\cal M}\cap
X$ and every $t\in I$.}
\vskip 3mm
As usual, $I$ is the unit interval $[0,1]$. 

The terminal homeomorphism $$h\colon\ {\cal L}\cap X\rightarrow
{\cal L}\cap X$$ of $H$, defined by $h(x):=H(x,1)$, of course
leaves ${\cal M}\cap X$ pointwise fixed. Such a homeomorphism
$h$ is called a {\it geometric monodromy of~${\cal L}\cap X$
relative to~${\cal M}\cap X$ above~$\omega_i$}. 

Another choice of loop~$\omega_i$ within the same homotopy
class~$\bar\omega_i$ in \hbox{$\hbox{\a
P}^1\setminus\bigcup_i\lambda_i$} and another choice of isotopy
$H$ above $\omega_i$ as in Lemma 1.2 would give a geometric
monodromy isotopic to~$h$ within~${\cal L}\cap X$ by an isotopy
leaving~${\cal M}\cap X$ pointwise fixed. Thus, the isotopy
class of $h$ in ${\cal L}\cap X$ relative to
${\cal M}\cap X$ is wholly determined by the homotopy class
$\bar\omega_i$ of $\omega_i$ in $\hbox{\a
P}^1\setminus\bigcup_i\lambda_i$. 

\vskip 1cm\goodbreak

\centerline{\f 2. Homological variation operators}\vskip 3mm

Fix an index $i$, and consider a geometric monodromy $h$
of~${\cal L}\cap X$ relative to~${\cal M}\cap X$ above
$\omega_i$. Denote by $S_q({\cal L}\cap X)$ the abelian group of
singular\break\noindent $q$-chains of ${\cal L}\cap X$ with
integer coefficients, and by $h_q\colon\ S_q({\cal L}\cap
X)\rightarrow S_q({\cal L}\cap X)$ the chain homomorphism
induced by $h$. Since $h$ leaves ${\cal M}\cap X$ pointwise
fixed ({\it cf.}~Lem\-ma~1.2), it is easy to see that for every
relative $q$-cycle $z$ of ${\cal L}\cap X$ modulo~${\cal M}\cap
X$, the variation by $h_q$ of $z$, that is the chain
$h_q(z)-z$, is an absolute $q$-cycle of ${\cal L}\cap X$ ({\it
cf.}~[C3,~Lem\-ma~4.6]). Moreover, one has the following
lemma.
\vskip 3mm

{\f Lemma 2.1} ({\it cf.}~[C3, Lemma 4.8]).~-- {\it The
correspondence}
$$\eqalign{
\hbox{var}_{i,q}\colon\ H_{q}({\cal L}\cap X,{\cal M}\cap
X)&\rightarrow H_{q}({\cal L}\cap X)\cr [z]_{{\cal L}\cap
X,{\cal M}\cap X} &\mapsto [h_q(z)-z]_{{\cal L}\cap X}\cr}$$
{\it gives a well-defined homomorphism which depends only on the
homotopy class~$\bar\omega_i$ of $\omega_i$ in~$\hbox{\a
P}^1\setminus\bigcup_i\lambda_i$.}\vskip 3mm

This means that another choice of loop $\omega_i$ within the
same homotopy class~$\bar\omega_i$ in~$\hbox{\a
P}^1\setminus\bigcup_i\lambda_i$ and another choice of monodromy
$h$ above $\omega_i$ as in Lemma~1.2 would give a homomorphism
equal to $\hbox{var}_{i,q}$.

Homomorphism $\hbox{var}_{i,q}$ is called a {\it homological
variation operator associated to $\bar\omega_i$}.

These operators were used by the first author in~[C3] to prove
the following result. 
\vskip 3mm
{\f Theorem 2.2} ({\it cf.}~[C3, Proposition 5.1]).~-- {\it If
$X$ is non-singular, then there is a natural isomorphism}
$$H_{d-1}({\cal L}\cap X)\Big/\sum_i \hbox{Im var}_{i,d-1}
\buildrel\sim\over\longrightarrow H_{d-1}(X).$$
{\it \indent In the special case of the complement of a projective
algebraic set, i.e., if $X:=\hbox{\a P}^n\setminus Z$, this can
be improved into an isomorphism}
$$H_{n+c-2}({\cal L}\cap X)\Big/\sum_i\hbox{Im var}_{i,n+c-2}
\buildrel\sim\over\longrightarrow
H_{n+c-2}(X),$$
{\it where $c$~is the least of the codimensions of the
irreducible components of~$Z$.}

\vskip 1cm

\centerline{\f 3. Homotopical variation operators on the
complements of}
\centerline{\f projective hypersurfaces with isolated
singularities}
\vskip 3mm

Throughout this section, we work under the following hypotheses.
\vskip 3mm

{\f Hypotheses 3.1}.~-- We assume that $Y=\hbox{\a P}^n$, with
$n\geq 3$, and $Z=H$, where~$H$ is a (closed) algebraic
hypersurface of $\hbox{\a P}^n$, with degree~$k$, having
only {\it isolated singularities}. Thus, $X=\hbox{\a
P}^n\setminus H$. We also assume that~${\cal S}$ is the
Whitney~stratifi\-cation the strata of which are:
$\hbox{\a P}^n\setminus H$, the singular part~$H_{\hbox{\bbb
sing}}$ of $H$, and the non-singular part $H\setminus
H_{\hbox{\bbb sing}}$ of $H$. Being transverse to~$\cal S$ then
means avoiding the singularities of~$H$ and being transverse to
the non-singular part of~$H$. Observe that
${\cal M}\cap X\neq\emptyset$. We fix a base
point~$*\in{\cal M}\cap X$.
\vskip 3mm

In~[CL], a $k$-fold (unramified) holomorphic covering\vskip -2mm
$$p\colon\ X'\to X$$\vskip -1mm\noindent
is constructed, where $X':=Y'\setminus Z'$ is a (pathwise) connected
quasi-projective variety in~$\hbox{\a P}^{n+1}$. In fact,
$X'$~is the global Milnor fibre of the cone of~$\hbox{\a
C}^{n+1}$ corresponding to~$H$. Moreover, it is shown that there
is a Whitney stratification~${\cal S}'$ of~$Y'$ preserving~$Z'$
and a pencil~${\cal P}'$ in~$\hbox{\a P}^{n+1}$ with axis~${\cal
M}'$ transverse to~${\cal S}'$ such that
$$\eqalign{
p^{-1}({\cal M}\cap X)&={\cal M}'\cap X'\qquad\hbox{and}\cr
p^{-1}({\cal P}(\mu)\cap X)&={\cal P}'(\mu)\cap X'\quad\hbox{for every }
\mu\in\hbox{\a P}^1,\cr}$$
the member ${\cal P}'(\mu)$ of ${\cal P}'$ with parameter $\mu$ being 
transverse to ${\cal S}'$ if and only if ${\cal P}(\mu)$ is transverse 
to ${\cal S}$. Recall that ${\cal L}={\cal P}(\lambda)$, and put 
${\cal L}':={\cal P}'(\lambda)$.
Then, for each~$i$, pencil~${\cal P}'$ gives
rise to a homological variation operator
$$\hbox{var}'_{i,n-1}\colon\ H_{n-1}({\cal L}'\cap X',{\cal
M}'\cap X')\to H_{n-1}({\cal L}'\cap X')$$
associated to~$\bar\omega_i$, defined as in section~2.

Given an index~$i$ and a base point $\bullet\in p^{-1}(*)$, one
can then consider the following diagram:
$$\matrix{
 &H_{n-1}({\cal L}'\cap X',{\cal M}'\cap X')
&\build\hbox to 12mm{\rightarrowfill}^{\hbox{var$'_{i,n-1}$}}
&H_{n-1}({\cal L}'\cap X')\cr
 &&&\cr
 &\chi\left\uparrow\vbox to 6mm{}\right.
&&\left\uparrow\vbox to 6mm{}\right.\hat\chi\cr
 &&&\cr
 &\pi_{n-1}({\cal L}'\cap X',{\cal M}'\cap X',\bullet)
&&\pi_{n-1}({\cal L}'\cap X',\bullet)\cr
 &&&\cr
 &\tilde\pi\left\downarrow\vbox to 6mm{}\right.
&&\left\downarrow\vbox to 6mm{}\right.\pi\cr
 &&&\cr
 &\pi_{n-1}({\cal L}\cap X,{\cal M}\cap X,*)
&&\pi_{n-1}({\cal L}\cap X,*),\cr
}\leqno(3.2)$$\smallskip\noindent
where $\tilde\pi$ and~$\pi$ are induced by~$p$ and
where
$\chi$ and~$\hat\chi$ are {\it Hurewicz homomorphisms}. Now, by a
general property of covering projections, $\tilde\pi$~is an
isomorphism ({\it cf.}~[Sp, Theorem~7.2.8]). Moreover,
homomorphism~$\hat\chi$ too is an isomorphism as a consequence
of the special fact that $X$~is the complement of a projective
hypersurface~$H$ with isolated singularities. Indeed,
${\cal L}\cap H$ is then a non-singular hypersurface of
${\cal L}\simeq\hbox{\a P}^{n-1}$, so that
$$\pi_1({\cal L}\cap X,*)\simeq \hbox{\a Z}/k\hbox{\a Z}\qquad
\hbox{and}\qquad
\pi_q({\cal L}\cap X,*)=0 \quad \hbox{for}\quad 2\leq q\leq
n-2$$
(this range may be empty) ({\it cf.}~[Li, Lemma 1.1]).
Knowing that ${\cal L}'\cap X'$ is pathwise
connected ({\it cf}.~[CL, Lemma~2.9]), these facts imply that
${\cal L}'\cap X'$ is
$(n-2)$-connected, and
$\hat\chi$ is then an isomorphism by the Hurewicz Isomorphism
Theorem.

Thus, for each $i$, and for every
$\bullet\in p^{-1}(*)$, there is a homomorphism
$${\cal V\hskip-0,8mm AR}_{i,n-1}\colon\
\pi_{n-1}({\cal L}\cap X,{\cal M}\cap X,*)\rightarrow
\pi_{n-1}({\cal L}\cap X,*)$$
defined by the composition
$$\pi\circ{\hat\chi}^{-1}\circ
\hbox{var}'_{i,n-1}\circ\chi\circ{\tilde\pi}^{-1}$$
in diagram~(3.2) ({\it cf}.~[CL, Section~5]).

One shows easily that homomorphism
${\cal V\hskip -0,8mm AR}_{i,n-1}$ does not in fact depend on the
choice of the base point $\bullet\in p^{-1}(*)$. 

Homomorphism ${\cal V\hskip -0,8mm AR}_{i,n-1}$ is called a {\it
homotopical variation operator associated to
$\bar\omega_i$}.

These operators were used in~[CL] to prove the following result.
\vskip 3mm

{\f Theorem 3.3} ({\it cf.}~[CL, Theorem 7.1]).~-- {\it Under
Hypotheses 3.1, there
is a natural isomorphism}
$$\pi_{n-1}({\cal L}\cap X,*)\Big/\sum_i \hbox{Im }{\cal V\hskip
-0,8mm
AR}_{i,n-1}\buildrel\sim\over\longrightarrow\pi_{n-1}(X,*).$$
\smallskip

As mentioned in the introduction, Theorem 3.3 is equivalent to
the projective version of [Li, Theorem~2.4] and also
provides a homotopy version of Theorem~2.2 in the special
case of complements of projective hypersurfaces with
isolated singularities.
\vskip 7mm

\centerline{\f 4. Generalized homotopical variation
operators}
\vskip 3mm

In this section, $X:=Y\setminus Z$ is again a (possibly
singular) quasi-projective variety as in Sections 1 and~2. We
assume further that ${\cal M}\cap X\not=\emptyset$ and we fix a
base point~$*$ in ${\cal M}\cap X$. Observe that the condition
${\cal M}\cap X\not=\emptyset$ is equivalent to $\dim X\geq 2$.
We also fix an index~$i$, and
consider a geometric monodromy $h$ of ${\cal L}\cap X$ relative
to ${\cal M}\cap X$ above $\omega_i$ ({\it cf.}~Section~1).

Let $q$ be an integer $\geq 1$. Since $h$~leaves ${\cal M}\cap X$
pointwise fixed, if
$f\in F^q({\cal L}\cap X,{\cal M}\cap X,*)$
({\it cf}.~Notation~0.1), then the map
$f\bot (h\circ f)$ defined on $I^q:=[0,1]^q$ by
$$f\bot (h\circ f)(x_1,\ldots,x_q):=\left\{\eqalign{
f(x_1,\ldots,x_{q-1},1-2x_q),\quad 0\leq x_q\leq\hbox{$1\over
2$},\cr h\circ f(x_1,\ldots,x_{q-1},2x_q-1),\quad \hbox{$1\over
2$}\leq x_q\leq 1,\cr}\right.$$
is in $F^q({\cal L}\cap X,*)$. Notice that the reversion
of~$f$ and its concatenation with $h\circ f$ are performed on
the variable transverse to the free face. This would in general
not make sense but here it does because $f$~and $h\circ f$ have
the same boundary.
\vskip 3mm

{\f Lemma 4.1.}~-- {\it The
correspondence}
$$\eqalign{
\hbox{VAR}_{i,q}\colon\ \pi_{q}({\cal L}\cap X,{\cal M}\cap
X,*)&\rightarrow
\pi_{q}({\cal L}\cap X,*)\cr
\langle f\rangle_{{\cal L}\cap X,{\cal M}\cap X,*} &\mapsto
\langle f\bot (h\circ f)\rangle_{{\cal L}\cap X,*}\cr}$$
{\it gives a well-defined map which depends only on the homotopy
class~$\bar\omega_i$ of $\omega_i$ in~$\hbox{\a
P}^1\setminus\bigcup_i\lambda_i$. If $q\geq 2$, it is a
homomorphism.}
\vskip 3mm

The independence assertion
follows from the remark we made just after Lemma~1.2. That the map is
well-defined is straightforward and one checks
easily that it is a homomorphism if $q\geq2$ (the sum of
homotopy cells being performed as in~[St, \S15]).

We shall call map $\hbox{VAR}_{i,q}$ a {\it
generalized homotopical variation operator associated
to~$\bar\omega_i$}. This terminology is justified by Theorem 5.1
below which asserts that, in the case where the homotopical
variation operators of~[CL] are defined ({\it cf.}~Section 3), the latter coincide with
our generalized operators.

We remark that if our operators are applied to absolute cells of
${\cal L}\cap X$ or if their result is considered as relative
cells of ${\cal L}\cap X$ modulo ${\cal M}\cap X$, then they act
as what can be called ordinary variations by monodromy. More
precisely:
\vskip3mm

{\f Observation 4.2.}~-- {\it Let\/}
$\hbox{incl}_q\colon\ \pi_{q}({\cal L}\cap X,*)
\rightarrow \pi_{q}({\cal L}\cap X,{\cal M}\cap X,*)$ {\it be the
natural map. Then,}\smallskip
\item{(i)} $\hbox{VAR}_{i,q}(\hbox{incl}_q(x))
               =-x+h_q(x)\ \hbox{for all\
        $x\in\pi_{q}({\cal L}\cap X,*)$};$\smallskip
\item{(ii)} if $q\geq 2$, $\ \hbox{incl}_q(\hbox{VAR}_{i,q}(y))
               =-y+\tilde h_q(y)\ \hbox{for all\
 $y\in\pi_{q}({\cal L}\cap X,{\cal M}\cap
X,*)$};$\smallskip\noindent {\it where $h_q$ and~$\tilde h_q$ are
the automorphisms of
$\pi_{q}({\cal L}\cap X,*)$ and $\pi_{q}({\cal L}\cap X,{\cal
M}\cap X,*)$ respectively induced by~$h$.}
\vskip3mm

The right-hand sides of the equalities are written additively
though the first group is not {\it a priori\/} commutative if
$q=1$ nor is the second one if $q=2$; the order of operations must
then be respected.

Observation~4.2 relies on the same reasons as those which allow the
sum of two homotopy cells to be performed indiscriminately on any
variable not transverse to the (possible) free face and which
make the sum commutative in high dimension. Its detailed
proof is left to the reader.

Notice that when $n=2$ and $X$~is the complement
$\hbox{\a P}^2\setminus C$ of a plane projective curve, ${\cal
M}\cap X$ is reduced to a single point and the observation above shows
that the epimorphism~(4), second row, of the introduction (the
existence of which will be justified by Lemma~4.8 below)
coincides with the map that van~Kampen's theorem asserts to be an
isomorphism.

Operator $\hbox{VAR}_{i,q}$ is linked to the homological
variation operator $\hbox{var}_{i,q}$ of Section~2 by Hurewicz
homomorphisms. This is stated in the next lemma.
\vskip 3mm

{\f Lemma 4.3.}~-- {\it The following diagram
is commutative:}
$$\matrix{
&H_{q}({\cal L}\cap X,{\cal M}\cap X) &\build\hbox to
12mm{\rightarrowfill}^{\hbox{var$_{i,q}$}} & H_{q}({\cal L}\cap
X)\cr &&&\cr &\chi\left\uparrow\vbox to 6mm{}\right. 
&&\left\uparrow\vbox to 6mm{}\right.\hat\chi\cr &&&\cr  
&\pi_{q}({\cal L}\cap X,{\cal M}\cap X,*)
&\build\hbox to 12mm{\rightarrowfill}^{\hbox{VAR$_{i,q}$}}
&\pi_{q}({\cal L}\cap X,*),\cr}$$
{\it where $\chi$ and $\hat\chi$ are Hurewicz homomorphisms.}
\vskip3mm

{\it Proof.} Since homotopy cells are defined on cubes, 
it is convenient to use {\it cubical\/} singular
homology theory ({\it cf.}~[HW], [M]), which is equivalent to
ordinary (simplicial) singular theory ({\it cf}.~[HW,
Section~8.4]). 
So, let us first introduce some
notation. For any pair of spaces
$(U,V)$, with
$V\subset U$, we shall denote by~$H_{q}^c(U,V)$ the
$q$-th cubical singular relative homology group of $(U,V)$, by
$[\cdot]^c_{U,V}$ the homology classes in this group, and by
$S_{q}^c(U,V)$ the group of $q$-dimensional cubical chains of
the pair~$(U,V)$. Given a (continuous) map
$g\colon\ (U,V)\rightarrow (U',V')$, we shall denote
by~$g_{q}\colon\ S_{q}^c(U,V)\rightarrow S_{q}^c(U',V')$ the
induced cubical chain homomorphism. A similar notation is used for
the absolute case.

Let~$f$ be a representative of an element of $\pi_q({\cal L}\cap
X,{\cal M}\cap X,*)$. We have
$$\eqalign{\chi(\langle f\rangle_{{\cal L}\cap
X,{\cal M}\cap X,*})&=[f_q(\iota)]_{{\cal L}\cap
X,{\cal M}\cap X}^c\cr
\hat\chi(\langle f\bot (h\circ f)\rangle_{{\cal L}\cap
X,*})&=[(f\bot (h\circ f))_q(\iota)]_{{\cal L}\cap X}^c\cr}$$
where $\iota\colon\ I^q\rightarrow I^q$ is the identity map ({\it
cf}.~[HW, 8.8.4]). 
Since
the expression for $\hbox{var}_{i,q}$ (given by Lemma~2.1) remains
valid in cubical theory (by [HW, 8.4.7 and the 
paragraph before~8.4.10]), it then suffices
to prove that the following equality holds in~$H_q^c({\cal
L}\cap X)$:
$$[h_q(f_q(\iota))-f_q(\iota)]_{{\cal L}\cap X}^c=
[(f\bot (h\circ f))_q(\iota)]_{{\cal L}\cap X}^c.\leqno (4.4)$$

For this purpose, consider the singular $q$-cubes
$\sigma_1,\sigma_2\colon\ I^q\rightarrow I^{q}$ in $I^{q}$
defined by
$$\eqalign{
\sigma_1(x_1,\ldots,x_{q})&:=(x_1,\ldots,x_{q-1},\hbox{$1-x_{q}\over
2$}),\cr
\sigma_2(x_1,\ldots,
x_{q})&:=(x_1,\ldots,x_{q-1},\hbox{$1+x_{q}\over 2$}).\cr}$$
Observe that $-\sigma_1+\sigma_2$ is a relative cycle of~$I^{q}$
modulo~$\dot I^{q}$.
\vskip 3mm

{\f Lemma 4.5.}~-- {\it The following equality holds in
$H_{q}^c(I^{q},\dot I^{q})$: 
$$[-\sigma_1+\sigma_2]^c_{I^{q},\dot I^{q}}=
[\iota]^c_{I^{q},\dot I^{q}}.$$}
\vskip 3mm

{\it Proof.} Let $\sigma\colon\ I^{q+1}\rightarrow I^{q}$ be the
singular
$(q+1)$-cube in $I^{q}$ defined by
$$\sigma(x_1,\ldots, x_{q+1}):=\left\{\eqalign{ (x_1,\ldots,
x_{q-1},\hbox{$2x_{q+1}+x_{q}-1\over 2$}),\cr (x_1,\ldots,
x_{q-1},\hbox{$1\over 2$}),\cr (x_1,\ldots,
x_{q-1},\hbox{$1-x_{q}+2x_{q+1}\over 2$}),\cr}\right.
\quad\eqalign{ &x_{q+1}\geq -\hbox{$x_{q}\over 2$} + 1,\cr
\hbox{$x_{q}\over 2$}\leq &x_{q+1}\leq -\hbox{$x_{q}\over 2$} +
1,\cr &x_{q+1}\leq \hbox{$x_{q}\over 2$}.\cr}$$  The boundary
operator
$$\partial\colon\ S_{q+1}^c(I^{q},\dot I^{q})\rightarrow
S_{q}^c(I^{q},\dot I^{q}),$$ applied to $\sigma$, satisfies
$$\partial \sigma=(-1)^{q+1}\sigma_1 -
(-1)^{q+1}\sigma_2-(-1)^{q}\iota.$$
Indeed, the face of~$\sigma$ of index~$(q,0)$ is degenerated and all
other non-mentioned faces are in~$\dot I^{q}$. The equality in the
statement of Lemma 4.5 follows.
\vskip 3mm

Now, since $f\bot (h\circ f)$~maps $(I^{q},\dot I^{q})$
into $({\cal L}\cap X,*)$, this lemma implies
$$[(f\bot (h\circ f))_{q}(\iota)]_{{\cal L}\cap X}^c =[(f\bot
(h\circ f))_{q}(-\sigma_1+\sigma_2)]_{{\cal L}\cap X}^c,$$ and
since $$(f\bot (h\circ f))\circ \sigma_1=f\quad\hbox{and}\quad
(f\bot (h\circ f))\circ \sigma_2=h\circ f,$$
one sees immediately
that
$$[(f\bot (h\circ f))_{q}(-\sigma_1+\sigma_2)]_{{\cal L}\cap X}^c
=[(h\circ f)_{q}(\iota)-f_{q}(\iota)]_{{\cal L}\cap X}^c.$$
This completes the proof of (4.4) and, consequently, the proof of
Lemma 4.3.
\vskip3mm

Operator $\hbox{VAR}_{i,q}$ also satisfies the following
equivariance property.
\vskip 3mm

{\f Lemma 4.6.}~-- {\it If $\gamma\in F^1({\cal M}\cap X,*)$ and
$f\in F^q({\cal L}\cap X,{\cal M}\cap X,*)$, then}
$$\hbox{VAR}_{i,q}(\langle \gamma\rangle_{{\cal M}\cap X,*}\cdot
\langle f\rangle_{{\cal L}\cap X,{\cal M}\cap X,*})=
\langle \gamma\rangle_{{\cal L}\cap X,*}\cdot
\hbox{VAR}_{i,q}(\langle f\rangle_{{\cal L}\cap X,{\cal M}\cap
X,*}),$$ {\it where $\cdot$ denotes equally the action of
$\pi_1({\cal M}\cap X,*)$ on 
$\pi_q({\cal L}\cap X,{\cal M}\cap X,*)$ or the action of 
$\pi_1({\cal L}\cap X,*)$ on $\pi_q({\cal L}\cap X,*)$.}
\vskip3mm

{\it Proof.} Let $\gamma^-$~be the inverse loop of~$\gamma$ and
$$K_{\gamma^-}\colon\ (I^q,\dot I^q)\times I\rightarrow ({\cal
L}\cap X, {\cal M}\cap X)$$
be a
$\gamma^-$-homotopy starting at $f$ (i.e.,
$K_{\gamma^-}(x,0)=f(x)$ for every
$x\in I^q$, $K_{\gamma^-}(x,t)=\gamma^-(t)$ for every
$x\in \dot I^q\setminus\{x_q=0\}$ and
every $t\in I$, and
$K_{\gamma^-}(x,t)\in {\cal M}\cap X$ for every $x\in
\{x_q=0\}$ and every $t\in I$). Denote by
$g$ the element of $F^q({\cal L}\cap X,{\cal M}\cap X,*)$ defined
by~$g(x):=K_{\gamma^-}(x,1)$. One has
$$\langle \gamma\rangle_{{\cal M}\cap X,*}\cdot
\langle f\rangle_{{\cal L}\cap X,{\cal M}\cap X,*}=
\langle g\rangle_{{\cal L}\cap X,{\cal M}\cap X,*}.\leqno (4.7)$$
Since $h$ leaves~${\cal M}\cap X$ pointwise fixed, the map
$${\cal K}_{\gamma^-}\colon\ (I^q,\dot I^q)\times I\rightarrow
({\cal L}\cap X, {\cal M}\cap X)$$ defined by
$${\cal K}_{\gamma^-}((x_1,\ldots,x_q),t):=\left\{\eqalign{
K_{\gamma^-}((x_1,\ldots,x_{q-1},1-2x_q),t), \quad 0\leq x_q\leq
\hbox{$1\over 2$},\cr  h\circ
K_{\gamma^-}((x_1,\ldots,x_{q-1},2x_q-1),t),\quad \hbox{$1\over
2$}\leq x_q\leq 1,\cr}\right.$$ is a $\gamma^-$-homotopy from
$f\bot (h\circ f)$ to $g\bot (h\circ g)$ such that
${\cal K}_{\gamma^-}(x,t)=\gamma^-(t)$ for every $(x,t)\in \dot
I^q\times I$. In other words,
$$\langle \gamma\rangle_{{\cal L}\cap X,*}\cdot
\hbox{VAR}_{i,q}(\langle f\rangle_{{\cal L}\cap X,{\cal M}\cap
X,*})=
\hbox{VAR}_{i,q}(\langle g\rangle_{{\cal L}\cap X,{\cal M}\cap
X,*}).$$ Lemma 4.6 then follows from (4.7).
\vskip3mm

Finally we prove a lemma which justifies the existence of the
epimorphisms (4) and~$(4')$ of the introduction. But this lemma
is valid for every $q\geq1$ and even if $X$~is singular.
\vskip3mm

{\f Lemma 4.8.}~-- {\it The image of operator\/}
$\hbox{VAR}_{i,q}$ {\it is contained in the kernel of the
natural map $\pi_q({\cal L}\cap X,*)\rightarrow\pi_q(X,*)$.}
\vskip 3mm

{\it Proof.} A representative of an element of the image of
$\hbox{VAR}_{i,q}$ is of the form $f\bot (h\circ f)$ with $f\in
F^q({\cal L}\cap X,{\cal M}\cap X,*)$. Let $H$ be an isotopy
giving rise to~$h$ as in Lemma~1.2. One defines a homotopy
$I^q\times I\rightarrow X$ from $f\bot (h\circ f)$ to the
constant map equal to $*$ by
$$\bigl((x_1,\ldots,x_q),t\bigr)\mapsto
\left\{\eqalign{
H\bigl(f(x_1,\ldots,x_{q-1},1-2x_q),2t\bigr),\cr
h\circ f(x_1,\ldots,x_{q-1},2x_q-1),\cr
h\circ f(x_1,\ldots,x_{q-1},1-2x_q),\cr
h\circ f(x_1,\ldots,x_{q-1},2t-1),\cr
h\circ f(x_1,\ldots,x_{q-1},2x_q-1),\cr}\right.
\quad\eqalign{
0\leq &t\leq\hbox{$1\over 2$},\cr
0\leq &t\leq\hbox{$1\over 2$},\cr
\hbox{$1\over 2$}\leq &t\leq 1,\cr
\hbox{$1\over 2$}\leq &t\leq 1,\cr
\hbox{$1\over 2$}\leq &t\leq 1,\cr}
\quad\eqalign{
0\leq &x_q\leq\hbox{$1\over 2$},\cr
\hbox{$1\over 2$}\leq &x_q\leq 1,\cr
0\leq &x_q\leq 1-t,\cr
1-t\leq &x_q\leq t,\cr
t\leq &x_q\leq 1.\cr}$$
By the first half of this homotopy, the lower part of the cell
undergoes the monodromy~$h$ while remaining attached to the
upper part; at the end of this process the two half cells become
opposite and the second half of the homotopy
collapses them together.

\vskip 3mm

{\f Remark.}~-- Since $\hbox{Im } H\subset X_{\partial K_i}$,
the proof above shows in fact that the image of
$\hbox{VAR}_{i,q}$  is contained in the kernel of the natural
map $\pi_q({\cal L}\cap X,*)\rightarrow\pi_q(X_{\partial
K_i},*)$.

\vskip 1cm

\centerline{\f 5. The link between ${\cal V\hskip -0,8mm
AR}_{i,n-1}$ and $\hbox{VAR}_{i,n-1}$ \ -- \  Main result}
\vskip 3mm

Throughout this section, we work under Hypotheses 3.1.
\vskip 3mm

{\f Theorem 5.1.}~-- {\it Under Hypotheses 3.1, the homotopical
variation operator ${\cal V\hskip -0,8mm AR}_{i,n-1}$ of
Ch\'eniot-Libgober which is then well-defined (cf.~Section 3)
coincides with the generalized homotopical variation operator\/}
$\hbox{VAR}_{i,n-1}$ {\it (defined in Section~4).}
\vskip 3mm
Before proving this theorem, observe that, together with
Theorem~3.3, it implies the following result.
\vskip 3mm

{\f Theorem 5.2.}~-- {\it Under Hypotheses 3.1, there is a
natural isomorphism}
$$\pi_{n-1}({\cal L}\cap X,*)\Big/\sum_i \hbox{Im
}\hbox{VAR}_{i,n-1}
\buildrel\sim\over\longrightarrow\pi_{n-1}(X,*).$$
\smallskip

Of course, Theorem 5.2 is equivalent to Theorem 3.3 and to
the projective version of [Li, Theorem~2.4]. 

We now turn to the proof of Theorem 5.1.
\vskip 3mm

{\it Proof of Theorem 5.1.} Consider the diagram obtained
from diagram~(3.2) by completing its lower row with the
homomorphism
$$\pi_{n-1}({\cal L}\cap X,{\cal M}\cap X,*)
\build\hbox to 12mm{\rightarrowfill}^{\hbox{VAR$_{i,n-1}$}}
\pi_{n-1}({\cal L}\cap X,*)$$
and its middle
row with the homomorphism
$$\pi_{n-1}({\cal L}'\cap X',{\cal M}'\cap X',\bullet)
\build\hbox to 12mm{\rightarrowfill}^{\hbox{VAR$'_{i,n-1}$}}
\pi_{n-1}({\cal L}'\cap X',\bullet)$$
defined from $\omega_i$ as
$\hbox{VAR}_{i,n-1}$ but with pencil~${\cal P}'$ and the
point~$\bullet$ instead of pencil~${\cal P}$ and the point~$*$.
We have to show that this new diagram is commutative.
But its lower square is indeed commutative since, given a
geometric monodromy~$h$ of~${\cal L}\cap X$ relative to~${\cal
M}\cap X$ above~$\omega_i$, there exists a geometric
monodromy~$h'$ of~${\cal L}'\cap X'$ relative to~${\cal M}'\cap
X'$ above~$\omega_i$ such that
$$p\circ h'=h\circ p$$
({\it cf}.~[CL, Remark~4.2]). As to the upper square, 
it commutes thanks to Lemma~4.3.

\vskip 1cm

\centerline{\f 6. A conjecture generalizing the van Kampen 
theorem to}
\centerline{\f non-singular quasi-projective varieties}\vskip3mm
In this section, we come back to the general hypotheses of
Section~1, so that $X:=Y\setminus Z$ is again a (possibly singular)
quasi-projective variety in~$\hbox{\a P}^n$ with $n\geq2$ as in
Sections~1, 2 and~4. We also assume that
${\cal M}\cap X\neq\emptyset$; as this condition is
equivalent to $\dim X\geq2$, it will be automatically fulfilled when $d\geq 2$.
We fix a base point~$*$ in
${\cal M}\cap X$.

The conjecture we made in the introduction can be specified
as follows.
\vskip 3mm
{\f Conjecture 6.1.}~-- {\it Under the hypotheses of Section~1 and
if $X$~is non-singular, there are  natural isomorphisms}
$$\eqalign{
\pi_{d-1}({\cal L}\cap X,*)\Big/\sum_i\hbox{Im VAR}_{i,d-1}
 &\buildrel\sim
  \over\longrightarrow\pi_{d-1}(X,*)\quad\hbox{\it if\/ $d\geq3$},\cr
\pi_1({\cal L}\cap X,*)\bigg/\overline{\bigcup_i\hbox{Im VAR}_{i,1}}
 &\buildrel\sim
  \over\longrightarrow\pi_1(X,*)\quad\hbox{\it if\/ $d=2$},\cr
}$$
{\it involving the generalized homotopical variation
operators\/} $\hbox{VAR}_{i,q}$ {\it defined in Section~4, with\/}
$\overline{\bigcup_i\hbox{Im VAR}_{i,1}}$ {\it denoting the normal
subgroup of $\pi_1({\cal L}\cap X,*)$ generated by\/}
$\bigcup_i\hbox{Im VAR}_{i,1}$.

{\it In the special case $Y=\hbox{\a P}^n$, so that
$X=\hbox{\a P}^n\setminus Z$, and provided $Z\not=\emptyset$, there
is a natural isomorphism}
$$\pi_{n+c-2}({\cal L}\cap X,*)\Big/\sum_i\hbox{Im
VAR}_{i,n+c-2}\buildrel\sim
\over\longrightarrow\pi_{n+c-2}(X,*)\quad
\hbox{\it if\/ $n+c\geq4$},$$
{\it  where $c$~is the
least of the codimensions of the irreducible components
of~$Z$ (notice that $n+c-2=d-1$ when $c=1$).}
\vskip 3mm

If proved, this conjecture would extend Theorem~5.2 (and hence
Theorem~3.3 and the projective version of~[Li, Theorem~2.4] reported
in the introduction as isomorphism~(2)) and would also extend the
classical Zariski-van~Kampen theorem on curves as remarked in
Section~4. It would give a complete homotopy analogue of Theorem~2.2
and thus would gather in a generalized form the Zariski-van~Kampen theorem
with a homotopical version of the Second Lefschetz Theorem.

We now give a first little approach to this conjecture.
\vskip 3mm

{\it First approach to Conjecture~6.1.} By Lemma~4.8, the subgroups
by which the quotients are taken are contained in the kernels of the
corresponding natural maps (which are epimorphisms by the Lefschetz 
Hyperplane Section Theorems, as mentioned in the introduction). 
The reverse inclusion which would lead to
the conclusion is much more difficult and not proved at the
moment.  We are simply giving below, via the following lemma, a
first little step in this direction.
\vskip 3mm

{\f Lemma 6.2.}~-- {\it If $X$~is non-singular and $d\geq2$, there is
a natural epimorphism
$$\pi_d(X_K,{\cal L}\cap
X,*)\rightarrow\pi_d(X,{\cal L}\cap X,*),$$
where $K$ is the union of the $K_i$ ({\it cf.}~Section 1). 

In the
special case where $Y=\hbox{\a P}^n$ and $Z\not=\emptyset$,
there is a natural epimorphism
$$\pi_{n+c-1}(X_K,{\cal L}\cap
X,*)\rightarrow\pi_{n+c-1}(X,{\cal L}\cap X,*),$$
where $c$~is as in Conjecture~6.1.}
\vskip 3mm

This lemma is a weak homotopical analogue of~[C3, Corollary~3.4]. It
shows that (with the same hypotheses of course)
the kernels of the natural maps
$$\eqalign{
\pi_{d-1}({\cal L}\cap X,*)\rightarrow\pi_{d-1}(X,*)\quad
 &\hbox{and}\quad\pi_{d-1}({\cal L}\cap X,*)
  \rightarrow\pi_{d-1}(X_K,*)\cr
\noalign{\noindent$\Bigl($resp.}
\pi_{n+c-2}({\cal L}\cap X,*)
\rightarrow\pi_{n+c-2}(X,*)\quad
 &\hbox{and}\quad\pi_{n+c-2}({\cal L}\cap X,*)
  \rightarrow\pi_{n+c-2}(X_K,*)\Bigr)\cr
}$$
are the same. So, with the remarks above, Conjecture 6.1 reduces
to the following one.
\vskip 3mm

{\f Conjecture 6.3.}~-- {\it Under the hypotheses of Conjecture
6.1, the kernel of the natural map $\pi_{d-1}({\cal L}\cap X,*)
\rightarrow\pi_{d-1}(X_K,*)$ is contained in\/}
$\sum_i\hbox{Im VAR}_{i,d-1}$ {\it if $d\geq3$ and in\/}
$\overline{\bigcup_i\hbox{Im VAR}_{i,1}}$ {\it if $d=2$. 

In the
special case where $Y=\hbox{\a P}^n$ and $Z\not=\emptyset$, the
kernel of the natural map
$\pi_{n+c-2}({\cal L}\cap X,*)\rightarrow\pi_{n+c-2}(X_K,*)$ is
contained in\/}
$\sum_i\hbox{Im VAR}_{i,n+c-2}$ {\it if $n+c\geq4$}.
\vskip 3mm

To complete this section, it remains to prove Lemma 6.2.
\vskip3mm

{\it Proof of Lemma 6.2.} By the homotopy sequence of the triple
$$(X,X_K,{\cal L}\cap X),$$
it suffices to prove that the pair $(X,X_K)$ is $d$-connected
(resp.~$(n+c-1)$-con\-nected).

We start by noticing that the pair $({\cal L}\cap X,{\cal M}\cap X)$
is $(d-2)$-connected (resp.~$(n+c-3)$-connected). This is shown by
applying the Lefschetz Hyperplane Section Theorem for non-singular
quasi-projective varieties (resp.~for complements) to the section of
${\cal L}\cap X$ by the hyperplane~$\cal M$ of~$\cal L$. To check
the required hypotheses and verify that the conclusion is the
announced one, we refer to the beginning of the proof of~[C3,
Corollary~3.4]. We just point out here that the hypothesis
$Z\neq\emptyset$ is crucial to ensure that the codimension of ${\cal
L}\cap Z$ in~$\cal L$ is still~$c$.

Thus, to show that
$(X,X_K)$ is $d$-connected (resp.~$(n+c-1)$-con\-nected), it is
enough to prove the following result which in fact holds even if
$X$~is singular.
\vskip 3mm

{\f Lemma 6.4.}~-- {\it For this lemma, $X$~may be singular. Let
$k$ be an integer~$\geq 0$. If
$({\cal L}\cap X,{\cal M}\cap X)$ is $k$-con\-nected, then
$(X,X_K)$ is $(k+2)$-connected.}
\vskip 3mm

This is a weak homotopy analogue of~[C3, Lemma~3.9]. In its proof,
the homology excision property is replaced by a much more
restrictive homotopy excision theorem, and the Eilenberg-Zilber
theorem and K\"unneth formula by a criterion on the connectivity of
the product of two pairs of spaces.
\vskip 3mm

{\it Proof of Lemma~6.4.} Let $\widetilde{\hbox{\a P}^n}$ be the blow
up of $\hbox{\a P}^n$ along ${\cal M}$, which is defined by
$$\widetilde {\hbox{\a P}^n}:=\{(x,\mu)\in \hbox{\a P}^n\times
\hbox{\a P}^1\mid x\in {\cal P}(\mu)\}.$$ It is a compact analytic
submanifold of
$\hbox{\a P}^n\times  \hbox{\a P}^1$ with dimension $n$. 

The projections of
$\hbox{\a P}^n\times\hbox{\a P}^1$ onto its two factors, when
restricted to
$\widetilde{\hbox{\a P}^n}$, give two proper analytic morphisms
$$\varphi\colon\ \widetilde{\hbox{\a P}^n}\rightarrow\hbox{\a
P}^n\quad
\hbox{and}\quad\pi\colon\ \widetilde{\hbox{\a
P}^n}\rightarrow\hbox{\a P}^1.$$ 

For any subsets $G\subset \hbox{\a P}^n$ and $E\subset \hbox{\a
P}^1$, write
$$\tilde G:=\varphi^{-1}(G)\quad \hbox{and}\quad \tilde G_E :=
\tilde G\cap \pi^{-1}(E).$$ One must not confuse $\tilde G_E$
with $\widetilde{G_E}$; we have
$$\widetilde{G_E}=\tilde G_E\cup \widetilde{(G\cap {\cal M})}
=\tilde G_E\cup(G\times\hbox{\a P}^1).$$

For simplicity, we also set $$L:={\cal L}\cap X\quad
\hbox{and}\quad M:={\cal M}\cap X.$$ 

 By stratifying suitably~$\tilde Y$ and then applying the First
Isotopy Theorem of Thom-Mather ({\it cf.}~[Th],~[Ma]) one obtains the
following lemma.
\vskip 3mm

{\f Lemma 6.5} ({\it cf.}~[C2, (11.1.5)]).~-- {\it The restriction
of $\pi$ to $\tilde X\setminus\bigcup_{i} \tilde X_{\lambda_i}$ is a
topological locally trivial fibration over
$\hbox{\a P}^1\setminus \bigcup_i \lambda_i$ with typical
fibre $\tilde X_\lambda$ homeomorphic to $L$. Moreover, this
bundle has $M\times (\hbox{\a P}^1\setminus \bigcup_i
\lambda_i)$ as a trivial subbundle of it.}
\vskip 3mm

The proof is now along lines very similar to~[La, 8.3]. Decom\-pose
$\hbox{\a P}^1$ into two closed semi-analytic hemispheres $D_+$ and $D_-$
such that: (i)~$D_+\cap D_-=S^1$, where
$S^1$ is a $1$-sphere; (ii) $K$~is contained in~$D_+$;
(iii)~$\lambda\in S^1$; (iv)~$D_-$ is contained in some
coordinate neighbourhood of the fibre bundle~$\tilde X\setminus
\bigcup_{i} \tilde X_{\lambda_i}$ for a trivialization which
preserves the subbundle $M\times (\hbox{\a P}^1\setminus \bigcup_i
\lambda_i)$. This choice of~$D_-$ implies that the pairs
$$(L\times D_-,L\times S^1\cup M\times D_-)\quad \hbox{and}\quad
(\tilde X_{D_-},\tilde X_{S^1}\cup\tilde M_{D_-})$$  are
homeomorphic.  Now, consider the following sequence of pairs of
spaces:
$$\eqalign{ (L\times D_-,L\times S^1\cup M\times D_-)&\simeq
(\tilde X_{D_-},\tilde X_{S^1}\cup\tilde
M_{D_-})\buildrel{exc}\over\hookrightarrow\cr &(\tilde X,\tilde
X_{D_+}\cup\tilde M)\hookleftarrow (\tilde X,\tilde
X_K\cup\tilde M)\buildrel{\varphi_\mid}\over\rightarrow
(X,X_K).\cr}$$

Since $(L,M)$ and $(D_-,S^1)$ are respectively $k$-connected and
$1$-connected,\break\noindent Exercise 9, p.~95 of [W] shows that
$(L\times D_-,L\times S^1\cup M\times D_-)$, and consequently
$(\tilde X_{D_-},\tilde X_{S^1}\cup\tilde M_{D_-})$, are
$(k+2)$-connected.

Next, for excision~$exc$, we might use the Homotopy Excision Theorem of
Blakers-Massey ({\it cf}.~[BM], [Sw, 6.21]) but we do not need its
full force. The elementary proof of~[La, (8.2.2)] which is about a
relative homeomorphism can be adapted to an excision and gives the
following result. Let $(A,B)\hookrightarrow(C,D)$ be an excision, that is
an inclusion of topological pairs such that $A\setminus B=C\setminus D$.
Suppose that $(A,B)$ is a relative CW-complex and that $A$ and~$D$ are
closed in~$C$. If $(A,B)$ is $m$-connected for some integer $m\geq0$,
then $(C,D)$ is also \hbox{$m$-connected}. Now, since
$(\tilde X_{D_-},\tilde X_{S^1}\cup\tilde M_{D_-})$ is
triangulable ({\it cf.}~[Lo]) and since $\tilde X_{D_-}$ and
$\tilde X_{D_+}\cup\tilde M$ are closed in~$\tilde X$, this result can be
applied to excision~$exc$ and shows that
$$(\tilde X,\tilde X_{D_+}\cup\tilde M)
\quad\hbox{is also $(k+2)$-connected}.\leqno (6.6)$$ 

Since $\tilde X\cap\pi^{-1} (D_+\setminus\bigcup_i\lambda_i)$ is
a fibre bundle ({\it cf.}~Lemma~6.5) and
$K\setminus\bigcup_i\lambda_i$ is a strong deformation retract
of $D_+\setminus\bigcup_i\lambda_i$, the First Homotopy Covering
Theorem [St, \S 11.3] shows that $\tilde X\cap\pi^{-1}
(K\setminus\bigcup_i\lambda_i)$ is a strong deformation retract
of 
$\tilde X\cap\pi^{-1} (D_+\setminus\bigcup_i\lambda_i)$.
Moreover, since the $\lambda_i$ are interior points of~$K$, the
deformation retraction may be in fact extended so that $\tilde
X_K$ is a strong deformation retract of $\tilde X_{D_+}$. As
furthermore $\tilde X_K$ is also a strong deformation retract of
$\tilde X_K\cup\tilde M_{D_+}$, one deduces that the pair
$(\tilde X_{D_+},\tilde X_K\cup\tilde M_{D_+})$ is
$q$-connected for all $q\geq 0$. By the theorem of Blakers-Massey
[Sw, Theorem 6.21], applied to the excision
$$(\tilde X_{D_+},\tilde X_K\cup\tilde M_{D_+}) \hookrightarrow
(\tilde X_{D_+}\cup\tilde M,\tilde X_K\cup\tilde M),$$
the same property
holds for
$(\tilde X_{D_+}\cup\tilde M,\tilde X_K\cup\tilde M)$. Then, the
homotopy sequence of the triple $$(\tilde X,\tilde X_{D_+}\cup
\tilde M,\tilde X_K\cup \tilde M),$$ together with (6.6),
implies that $(\tilde X,\tilde X_K\cup\tilde M)$ is
$(k+2)$-connected. 

Now, the same holds for $(X,X_K)$ since $$\varphi_\mid\colon\
(\tilde X,\tilde X_K\cup\tilde M)\rightarrow (X,X_K)$$ is a
relative homeomorphism ({\it cf.}~[La, (8.2.2)]). This completes
the proof of Lemma 6.4 and, consequently, the proof of Lemma
6.2.
\vskip 3mm

Conjecture~6.3, to which Conjecture~6.1 is thus reduced, remains of
course the hard part of the work.

The generalized homotopical variation operators introduced here 
are also certainly a first step for further generalizations of van Kampen's
theorem to {\it singular} quasi-projective varieties, the influence of 
the singularities being measured by the (local or global) rectified 
homotopical depth ({\it cf.}~[G], [HL1,2], [GM1,2], [E1,2]).

\vskip 8mm

{\it Acknowledgements.} The second author would like to thank
the ABDUS SALAM International Centre for Theoretical Physics
(ICTP) in Trieste for kind hospitality and support.

\vskip 8mm\goodbreak

\centerline{\f References}\vskip 3mm
\frenchspacing
\item{[AF]}  {\f A.~Andreotti, T.~Frankel},  ``The second
Lefschetz theorem on hyperplane sections'', {\it Global
analysis, papers in honor of K.~Kodaira\/}
(eds D.~C.~Spencer and S.~Iyanaga, University of Tokyo Press,
Tokyo; Princeton University Press, Princeton, NJ, 1969),
1--20.\smallskip
\item{[BM]} {\f A.L. Blakers, W.S. Massey}, ``The homotopy groups of a triad,
II'', {\it Ann.~of Math.}~(2)~{\bf 55} (1952)  192--201.\smallskip\goodbreak
\item{[C1]} {\f D. Ch\'eniot}, ``Une d\'emonstration du th\'eor\`eme de
Zariski sur les sections hyperplanes d'une hypersurface projective et du
th\'eor\`eme de van Kampen sur le groupe fondamental du compl\'ementaire
d'une courbe projective plane'',   {\it Compositio Math.}~{\bf 27} (1973)  
141--158.\smallskip
\item{[C2]} {\f D. Ch\'eniot}, ``Topologie du
compl\'emen\-taire  d'un ensemble alg\'ebrique projectif'',  
{\it Ensei\-gn. Math.}~(2) {\bf 37} (1991)  
293--402.\smallskip
\item{[C3]} {\f D. Ch\'eniot}, ``Vanishing cycles in a pencil
of hyperplane sections of a non singular quasi-projective
variety'',   {\it Proc.~London Math.~Soc.}~(3) {\bf 72} (1996)
515--544.
\smallskip
\item{[CL]} {\f D. Ch\'eniot, A. Libgober}, ``Zariski-van Kampen theorem for
higher homotopy groups'', LATP-Preprint {\bf 01-25}, UMR CNRS 6632,
Universit\'e de Provence, Marseille, 2001,
arXiv:math.AG/0203019.\smallskip
\item{[E1]} {\f C. Eyral}, ``Profondeur homotopique et conjecture de
Grothendieck'', {\it Ann.~Sci.~\'Ecole Norm.~Sup.} (4) {\bf 33} 
(2000) 823--836.
\smallskip
\item{[E2]} {\f C. Eyral}, ``Tomographie des vari\'et\'es singuli\`eres
et th\'eor\`emes de Lefschetz'', {\it Proc.~London Math.~Soc.} (3) {\bf 83} 
(2001) 141--175.
\smallskip
\item{[E3]} {\f C. Eyral}, ``Topology of quasi-projective
varieties and Lefschetz theory'', ICTP-Preprint {\bf
IC/2001/159}, The Abdus Salam International Centre for
Theoretical Physics, Trieste, 2001.
\smallskip
\item{[GM1]}  {\f M. Goresky, R. MacPherson},  ``Stratified
Morse theory'', {\it Singularities, Part 1 (Arcata, CA, 1981)},
Proceedings of Symposia in Pure Mathematics {\bf 40} (American
Mathematical Society, Providence, RI, 1983) 517--533.\smallskip
\item{[GM2]} {\f M. Goresky, R. MacPherson}, {\it Stratified
Morse theory\/} (Springer-Verlag, New-York, 1988).\smallskip
\item{[G]} {\f A. Grothendieck}, {\it Cohomologie locale des faisceaux
coh\'erents et th\'eor\`e\-mes de Lefschetz locaux et globaux (SGA2)\/} 
(Adv.~Stud.~Pure Math.~2, North-Holland, Amsterdam, 1968).\smallskip
\item{[HL1]} {\f H.A. Hamm, L\^e D.T.}, ``Lefschetz theorems on
quasi-projective varieties'',  {\it Bull. Soc.~Math. France\/}
{\bf 113} (1985) 123--142.\smallskip
\item{[HL2]} {\f H.A. Hamm, L\^e D.T.}, ``Rectified homotopical
depth and Grothendieck conjectures'',  {\it The Grothendieck
Festschrift}, vol.~II (Birkh\"auser, Boston, 1991)
311--351.\smallskip
\item{[HW]} {\f P.J. Hilton, S. Wylie}, {\it Homology theory - An introduction
to algebraic topology\/} (Cambridge University Press, Cambridge,
1965).\smallskip
\item{[La]}  {\f K. Lamotke},  ``The topology of complex projective varieties
after S. Lefschetz'', {\it Topology\/}~{\bf 20} (1981)
15--51.\smallskip
\item{[LT]}  {\f L\^e D.T., B. Teissier},   ``Cycles \'evanescents, sections
planes et conditions de Whitney II'', {\it Singularities, Part 2 (Arcata, CA,
1981)}, Proceedings of Symposia in Pure Mathematics {\bf 40} (American
Mathematical Society, Providence, RI, 1983)  65--103.\smallskip 
\item{[Li]}  {\f A. Libgober},  ``Homotopy groups of the complements to singular
hypersurfaces,~II'', {\it Ann.~of Math.}~(2) {\bf 139} (1994) 117--144.\smallskip
\item{[Lo]}  {\f S. \L ojasiewicz}, ``Triangulation of
semi-analytic sets'', {\it Annali Scu.~Norm. Sup.~Pisa\/} {\bf
18} (1964) 449--474.\smallskip
\item{[M]} {\f W.S. Massey}, {\it A basic course in algebraic
topology\/} (Graduate Texts in Mathematics {\bf 127},
Springer-Verlag, New-York, 1991).\smallskip
\item{[Ma]}  {\f J. Mather}, {\it Notes on topological
stability\/} (mimeographed notes, Harvard University,
1970).\smallskip
\item{[Sp]} {\f E.H. Spanier}, {\it Algebraic topology\/}
(Reprint of the 1966 original, Springer-Verlag,  New-York,
1989).\smallskip 
\item{[St]} {\f N. Steenrod}, {\it The topology of fibre
bundles\/} (Princeton University Press, Princeton,
1951).\smallskip 
\item{[Sw]} {\f R.M. Switzer},  {\it Algebraic topology - Homotopy and
homology\/}  (Springer-Verlag, New-York, 1975).\smallskip  
\item{[Th]} {\f R. Thom},  ``Ensembles et morphismes stratifi\'es'',
{\it Bull.~Amer.~Math.~Soc.}~{\bf 75} (1969)
240--284.\smallskip  
\item{[vK]}  {\f E.R. van Kampen},  ``On the fundamental group of an algebraic
curve'',  {\it Amer.~J. Math.}~{\bf 55} (1933) 255--260.\smallskip
\item{[W]}  {\f G.W. Whitehead}, {\it Elements of homotopy
theory\/} (Graduate Texts in Mathematics {\bf 61},
Springer-Verlag, New-York, 1978).\smallskip
\item{[Wh]}  {\f H. Whitney},  ``Tangents to an analytic
variety'', {\it Ann.~of Math.}~(2) {\bf 81} (1965) 496--549.
\smallskip
\item{[Za]}  {\f O.~Zariski},  ``On the problem of existence
of algebraic functions of two variables possessing a given
branch curve'', {\it Amer. J. Math.} {\bf 51} (1929)
305--328 .
\smallskip

\bye